\RequirePackage{fix-cm}
\documentclass[smallextended]{svjour3}       
\smartqed  
\usepackage{graphicx}

\usepackage{amsmath}
\usepackage{amsfonts}
\usepackage{amssymb}
\usepackage{algpseudocode}
\usepackage[ruled]{algorithm}

\newtheorem{assumption}{Assumption}

\DeclareMathOperator*{\argmin}{\arg\min}

\begin{document}

\title{A clustering heuristic to improve a derivative-free algorithm for nonsmooth optimization
}


\author{Manlio Gaudioso         \and
        Giampaolo Liuzzi \and
        Stefano Lucidi 
}


\institute{M. Gaudioso \at
              Universit\'a della Calabria, Dipartimento di Ingegneria Informatica, Modellistica, Elettronica e Sistemistica, 87030 Rende (CS), Italy. \\
              \email{manlio.gaudioso@unical.it}           
           \and
           G. Liuzzi \at
              ``Sapienza'' Universit\`a  di Roma, Dipartimento di Ingegneria Informatica Automatica e Gestionale, Via Ariosto 25, 00185 Rome, Italy.\\
              \email{liuzzi@diag.uniroma1.it}
           \and
			S. Lucidi \at
			``Sapienza'' Universit\`a  di Roma, Dipartimento di Ingegneria Informatica Automatica e Gestionale, Via Ariosto 25, 00185 Rome, Italy.\\
			\email{lucidi@diag.uniroma1.it}
}

\date{Received: date / Accepted: date}

\maketitle

\begin{abstract}
In this paper we propose an heuristic to improve the performances of the recently proposed derivative-free method for nonsmooth optimization CS-DFN. The heuristic is based on a clustering-type technique to compute a direction { which relies on an estimate of Clarke's generalized gradient} of the objective function. As such, this direction (as it is shown by the numerical experiments) is a good descent direction for the objective function. We report some numerical results and comparison with the original CS-DFN method to show the utility of the proposed improvement on a set of well-known test problems. 	
\end{abstract}

\section{Introduction}
We consider the following unconstrained minimization problem
\begin{equation}\label{prob}
	\min_{x\in\Re^n}\ f(x).
\end{equation}
We assume that the objective function $f$, (though nonsmooth) is
Lipschitz continuous and that first-order information is unavailable or impractical to obtain. {We require the following assumption.}
\begin{assumption}
	The function $f(x)$ is coercive, i.e. every level set is compact.
\end{assumption}
Useless to say that there is plenty of problems with the above features, especially coming from the engineering context. In the literature, many approaches have been proposed to tackle the nonsmooth problem (\ref{prob}) { in the derivative-free framework}. They can be roughly subdivided into two main classes: direct-type algorithms and model-based algorithms. 
\begin{itemize}
	\item[-] {\em Direct-type methods}. The algorithms belonging to this class make use of suitable sampling of the objective function. They occasionally can heuristically use modeling techniques, but the convergence theory hinges on the sampling technique. In this class of methods, we cite the mesh adaptive direct search algorithm implemented in the software package NOMAD \cite{nomad4,nomad3}, the linesearch derivative-free algorithm CS-DFN proposed in \cite{Fasano:14} { and the discrete gradient method \cite{bagirov:08}.}
	\item[-] {\em Model-based methods}. This class comprises all those algorithms whose convergence is based on the strategy used to build the approximating models. Within this class we can surely cite the recent trust-region derivative-free method proposed in \cite{liuzzi:19}.   
\end{itemize}

In the relatively recent paper \cite{Fasano:14}, a method for optimization of nonsmooth black-box problems has been proposed, namely CS-DFN. CS-DFN is able to solve problems more general than problem (\ref{prob}) above since it can  handle also nonlinear and bound constraints. It is based on a penalization approach, namely the nonlinear constraints are penalized by an exact penalization mechanism whereas (possible) bound constraints on the variables are handled explicitly.

In this paper, we propose an improvement of CS-DFN by incorporating into its main algorithmic scheme a clustering heuristic to compute efficient search directions. {Starting from an approximation of the directional derivatives along a certain set of directions, we construct a polyhedral approximation of the subdifferential which in turn is used to calculate a search direction in the steepest descent fashion. Along such direction we implement}  a linesearch procedure with extrapolation just like the one adopted by CS-DFN to explore its directions.

To asses the potentialities of the proposed improvement, we carry out an experimentation and comparison of CS-DFN with and without the proposed heuristic. The results, in our opinion, clearly show the advantages of the improved method over the original one.

The paper is organized as follows. In section \ref{sec:directions} we extend to a nonsmooth setting the steepest descent direction and a kind of Newton-type directions. In section \ref{sec:heuristic} we propose an heuristic to compute possibly efficient directions in a derivative-free context. Section \ref{sec:improvedalgo} we describe an improved version of the CS-DFN algorithm which is obtained by suitably employing the improved directions just described. In section \ref{sec:results} we report the results of a numerical comparison between CS-DFN and the proposed improved version on a set of well-known test problems. Finally, section \ref{sec:conclusion} is devoted to some discussion and conclusions.

\subsection{Definitions and notations}

\begin{definition}\label{def:Clarke}
	Given a point $x\in\Re^n$ and a direction $d\in\Re^n$, the Clarke directional derivative of $f$ at $x$ along $d$  is defined as \cite{clarkebook:90}
	\[
	f^\circ(x;d) = \limsup_{y\to x,t\downarrow 0}\frac{f(y+td) - f(y)}{t}.
	\]
	{Moreover the Clarke generalized gradient (or subdifferential) $\partial_C f(x)$ is defined as
		\[
		\partial_C f(x) = conv \{s :s \in \Re^n , \nabla f({x}_k) \rightarrow s, ~ x_k \rightarrow {x}, ~ x_k \not \in \Omega_f\},
		\] 
		$\Omega_f$ being the set (of zero measure) where $f$ is not differentiable.}
\end{definition}

{The following property holds:
	\begin{equation}
		f^\circ(x;d)= \max_{s \in \partial_C f(x)}s^{\top}d
\end{equation}}

\begin{definition}\label{def:stationarity}
	A point $x^*\in\Re^n$ is Clarke stationary for Problem (\ref{prob}) when $f^\circ(x^*,d) \geq 0$, for all $d\in\Re^n$.
\end{definition}

{In the following, we denote by $e_i$, $i=1,\dots,n$, the $i$-th column of the canonical basis in $\Re^n$ and by $e$ a vector of all ones of appropriate dimensions.}

\section{Descent type directions}\label{sec:directions}
In the context of nonsmooth optimization, efficient search directions can be computed by using the information provided by {the subdifferential} of the objective function. In the following subsections, we describe how such directions can be obtained.

\subsection{Steepest descent direction  $g_k^S$}
In this subsection we recall a {classic approach \cite{makela:92} to compute a generalization to nonsmooth functions of the steepest descent direction for continuously differentiable functions.} 

\noindent
Let us consider the vector which minimizes the following ``first order-type'' model of the objective function.
\begin{equation}\label{I-model}
	\begin{array}{l}
		\hbox{min}\ \ f(x_k)+f^\circ(x_k; d)+\frac{1}{2}\|d\|^2\\[0.4em]
	\end{array}
\end{equation}
\par\medskip\noindent
Note that, in the case of continuously differentiable functions, we have that $f^\circ(x_k; d)=\nabla f(x_k)^Td_k$ and that the solution of Problem (\ref{I-model}) is given by $d^*=-\nabla f(x_k)$.
\par\medskip\noindent
For  nonsmmoth functions, standard results \cite{makela:92} lead to the following proposition.
\par\medskip\noindent
\begin{proposition}
	Let $d^S$ be the solution of Problem (\ref{I-model}). Then
	\begin{description}
		\item[i)] the vector $d^*$ is given by
		\begin{equation}\nonumber
			d^S=-g_k^S
		\end{equation}
		where
		\begin{equation}\label{eq:gks}
			\begin{array}{l}
				g_k^S=\quad \hbox{argmin}\ \|\xi\|^2\\[0.4em]
				\qquad\qquad  s.t.\  \xi\in\partial f(x_k)
			\end{array}
		\end{equation}
		\item[ii)]the vector $d^S$ satisfies $f^\circ(x_k; -g_k^S)=-\|g_k^S\|^2$.
		\item[iii)] for any $\gamma\in (0,1)$ a $\bar \alpha$ exists such that
		$$f(x_k-\alpha g_k^S)\le f(x_k) - \alpha \gamma \|g_k^S\|^2$$
		with $\alpha\in (0, \bar \alpha]$
	\end{description}
\end{proposition}
\par\medskip\noindent

{
	The above $d_k^S$ direction is a first-order direction which (closely) resembles  the steepest-descent direction for continuously differentiable case.
}

\subsection{Newton-type direction  $d_k^N$}\par\medskip\noindent
In the nonsmooth case, obtaining a Newton-type direction is much more involved than in the differentiable case. In the latter case {it suffices to pre-multiply} the anti-gradient by the Hessian of the objective function. In the nonsmooth case {instead of simply  pre-multiplying direction $g_k^S$ by any positive definite matrix, we resort to  minimizing  the following ``second order-type'' model.}
\begin{equation}\label{II-model}
	\hbox{min}\ \ f(x_k)+f^\circ(x_k; d)+\frac{1}{2}d^TB_kd
\end{equation}
where $B_k$ is a positive definite matrix. Let us call the solution of problem (\ref{II-model}) $d_k^N$.
\par
\medskip
\noindent
For problem (\ref{II-model}) the following proposition can be proved.
\par\medskip\noindent
\begin{proposition}
	Let $d_k^N$ be the solution of Problem (\ref{II-model}). Then
	\begin{description}
		\item[i)] the vector $d_k^N$ is given by
		\begin{equation}\nonumber
			d_k^N=-B_k^{-1}g_k^N
		\end{equation}
		where
		\begin{equation}\label{sub-prob-new}
			\begin{array}{l}
				g_k^N=\quad \hbox{argmin}\ \xi^TB_k^{-1}\xi\\[0.4em]
				\qquad\qquad  s.t. \ \xi\in\partial f(x_k)
			\end{array}
		\end{equation}
		\item[ii)]the vector $d_k^N$ satisfies\qquad $f^\circ(x_k; d_k^N)=-(g_k^N)^TB_k^{-1}g_k^N=(g_k^N)^Td_k^N$.
		\item[iii)] for any $\gamma\in (0,1)$ a $\bar \alpha$ exists such that
		$$f(x_k-\alpha B_k^{-1}g_k^N)\le f(x_k) - \alpha \gamma (g_k^N)^TB_k^{-1}g_k^N$$
		with $\alpha\in (0, \bar \alpha]$
	\end{description}
\end{proposition}
\par\noindent
{\bf Proof}. By repeating the similar arguments of proof of Theorem 5.2.8 in \cite{makela:92} we have that function $\phi(d)=f^\circ(x_k; d)+\frac{1}{2}d^TB_kd$ is strictly convex. Therefore Problem  (\ref{II-model}) has a  unique minimizer $d^*$ such that:
\begin{eqnarray}\label{ott-1}
	&&0\in \partial f^\circ(x_k; d^*)+B_kd^*.
\end{eqnarray}
Recalling Lemma 5.2.7 of \cite{makela:92} we have:
\begin{eqnarray}\label{ott-2}
	&&\partial f^\circ(x_k; d^*)\subseteq \biggl\{\xi\in \partial f(x_k):\xi^Td^*= f^\circ(x_k; d^*)\biggr\}
\end{eqnarray}
The relations (\ref{ott-1}) and (\ref{ott-2}) imply that a vector $g_k^N$ exists such that:
\begin{eqnarray*}
	&&g_k^N=-B_kd^*,\\
	&&(g_k^N)^Td^*=  f^\circ(x_k; d^*).
\end{eqnarray*}
and, hence,
\begin{eqnarray}\label{ott-3}
	&&-(g_k^N)^TB_k^{-1}g_k^N=  f^\circ(x_k; -B_k^{-1}g_k^N)
\end{eqnarray}
which proves point ii) by setting $d_k^N=d^*$.
\par\smallskip\noindent
Now the definition of $ f^\circ(x_k; -B_k^{-1}g_k^N)$ and (\ref{ott-3}) give:
\begin{eqnarray*}
	&& f^\circ(x_k; -B_k^{-1}g_k^N)= \max_{\xi\in \partial f(x_k)}  \xi^T( -B_k^{-1}g_k^N) =-(g_k^N)^TB_k^{-1}g_k^N.
\end{eqnarray*}
\begin{eqnarray}\label{ott-4}
	&&(g_k^N)^TB_k^{-1}g_k^N\le   (g_k^N)^TB_k^{-1}\xi,\qquad\hbox{for all}\qquad \xi\in \partial f(x_k),
\end{eqnarray}
which implies
\begin{eqnarray}\label{ott-5}
	&&(g_k^N)^TB_k^{-1}(\xi-g_k^N)\ge  0,\qquad\hbox{for all}\qquad \xi\in \partial f(x_k),
\end{eqnarray}
Therefore, (\ref{ott-5}) shows that the vector $g_k^N$ is the unique solution of Problem \ref{sub-prob-new}.
\par\smallskip\noindent
Finally point iii) again follows from definition of $ f^\circ(x_k; -B_k^{-1}g_k^N)$ and (\ref{ott-3}).\hfill$\triangleleft$

\section{An heuristic approach to define efficient directions}\label{sec:heuristic}

At the base  of the proposed heuristics is the hypothesis that nonsmoothness of the objective function is due to its {\it  finite $\max$} structure. Such hypothesis appear realistic as a wide range of nonsmooth optimization problems, coming from practical applications, are of the $\min\max$ type.  Drawing inspiration from the paper \cite{lukvle:98} (see also \cite{frangioni11}), given points $y_j\in\Re^n$, $j= \{1,2,\dots,p\}$, sufficiently close to $x$, the (possibly) non-convex and non-smooth function $f(x)$ {is} approximated by using the following piece-wise quadratic model function,
\[
f^\Box(x) = \max_{j=1,\dots,p}\{q_j(x)\}
\]
with
\[
q_j(x) = f(y_j) + g_j^\top (x-y_j) + \frac{1}{2}(x-y_j)^\top H_j(x-y_j) 
\]
where $g_j\in\partial f(y_j)$ and $H_j = H(y_j)$, $j=1,\dots,p$. 
We remark that, while we assume that {the model  structure of $f$ is  a  } $\max$ of a finite number of functions, the number $p$ of such functions is unknown and has to be estimated via a trial--and--error calculation process.

We can write, 
\[
\partial f^\Box(x) = \partial\max_{j=1,\dots,p}\{q_j(x)\} \subseteq conv\{g_j + H_j(x-y_j), j=1,\dots,p\} = C(x).
\]
Furthermore, by assuming that $f(x)\approx f^\Box(x)$, we have
\begin{equation}\label{clarkederiv_approx}
	f^\circ(x;d) = \max_{s\in\partial f(x)}d^\top s \approx
	\max_{s\in\partial f^\Box(x)} d^\top s \leq \max_{s\in  C(x)} d^\top s = d^\top (g_{\bar\imath} + H_{\bar\imath}(x-y_{\bar\imath})). 
\end{equation}
\par\smallskip
In the actual case, $C(x)$ is the convex hull of a given number of generator vectors $v_j$, $j=1,\dots,p$. We can try and estimate those generators by using the quantities computed by the algorithm. 

More in particular,  let $x_k$ be the current iterate of the algorithm, $d_i\in\Re^n$ and $\alpha_i > 0$, $i=1,\dots,r$, the directions sampled by the algorithm along with their respective stepsizes, and define
\begin{equation}
	s_i = \frac{f(x_k+\alpha_id_i) - f(x_k)}{\alpha_i} \approx f^\circ(x_k;d_i). \label{incremental}
\end{equation} 
By using (\ref{clarkederiv_approx}), for $i=1,\dots,r$,
\[
f^\circ(x_k;d_i) \approx d_i^\top v_{j_i},\quad j_i\in\{1,2,\dots,p\}.
\]
It is then possible to compute estimates of the generators $v_j$, $j=1,\dots,p$, by solving the problem
\begin{equation}
	\min_{\hat v_1,\dots,\hat v_p}\sum_{i=1}^r \min_{j=1,\dots,p}\{(d_i^\top \hat v_j-s_i)^2\}. \label{minmin}
\end{equation}
The above problem 
{is a hard, nonsmooth nonconvex problem of the clustering type. It can be put however in  DC (Difference of Convex) form as in \cite{khalaf17}. Since it has to be solved many times during the proposed algorithm, we prefer to resort, in our implementation, to a greedy heuristic of the $k$-means-type  \cite{kmeans2,kmeans1,kmeans3}}. 

\begin{algorithm}[H]
	\caption{$k$-means-type Algorithm ($r,p,G$)}\label{kmeansAlgo}
	
	\begin{algorithmic}[1]
		\medskip
		
		\State   {\bf Data.} $r,p\in\mathbb{N}$, $G=\{(d_i,s_i), {i=1,\dots,r}\}$, $d_i\in\Re^n$, $s_i\in\Re$,  $h_{\max} > 0$.
		
		\medskip
		
		\State Set $I = \{1,\dots,r\}$, $\hat v_j^{(0)}=0$, $j=1,\dots,p$.
		
		\For{$h=0,\dots,h_{\max}$}
		
		\For{$j=1,\dots,p$} 
		\State Compute
		\begin{eqnarray*}
			I_j^{(h)} & = & \{i\in I: (d_i^\top \hat v_j^{(h)} - s_i)^2 \leq (d_i^\top \hat v_\ell^{(h)} - s_i)^2,\ \ell\neq j\}\\
			\hat v_j^{(h+1)} & = & \displaystyle\argmin_{v\in\Re^n} \sum_{i\in I_j^{(h)}} (v^\top d_i - s_i)^2 \\
			\phi_j^{(h+1)} & = & \displaystyle\sum_{i\in I_j^{(h)}} ((\hat v_j^{(h+1)})^\top d_i - s_i)^2 
		\end{eqnarray*} 
		\EndFor
		\EndFor
		
		\State {\bf Return} $\hat v_j^{(h_{\max}+1)}$ and $\phi_j^{(h_{\max}+1)}$ for $j=1,\dots,p$.

		\par\medskip\noindent
		
	\end{algorithmic}
\end{algorithm}

Then, we can compute an estimate of direction $d_k^N$ by solving problem (\ref{eq:gks}) (or (\ref{sub-prob-new})) where $\partial f(x_k)$ (or $\partial_\epsilon f(x_k)$) is approximated by $conv(\hat v_i,\dots,\hat v_p)$. More precisely, we define the following algorithm that computes a search direction.

\begin{algorithm}[H]
	\caption{Direction computation ($r,\tilde\alpha,y,B,G$)}\label{DircoAlgo}
	
	\begin{algorithmic}[1]
		\medskip
		
		\State {\bf Data}. $r\in\mathbb{N}$,  $G=\{(d_i,s_i), {i=1,\dots,r}\}$, $d_i\in\Re^n$, $s_i\in\Re$,  $B\in{\cal S}^{n\times n}$, $B\succ 0$.
		
		\State Set $\hat\alpha = 0$, $\hat d^N = 0$
		\For{$p=2,\dots,\min\{r,n\}$}
		
		\State Compute $\hat v_1,\dots, \hat v_p$ and $\phi_1,\dots,\phi_p$ by $k$-means-type Algorithm ($r,p,G$), i.e. Algorithm \ref{kmeansAlgo}.
		
		\If{$\sum_{i=1}^p\phi_i < \epsilon$}
		
		\State Compute $\hat d^N$ by
		
		\[
		\hat d^N = \argmin_{\xi\in conv(\hat v_1,\dots,\hat v_p)} \xi^\top B^{-1}\xi.
		\]
		
		
		
		\EndIf

		\EndFor
		\par\medskip\noindent
		
	\end{algorithmic}
\end{algorithm}

In the following we give an example of how the heuristic works.
\begin{example}
	{
		Consider the (convex) nonsmooth function $maxl$ \cite{lukvle:98}, defined as
		\[ f(x)= \max_{1 \leq i \leq n} |x_i|.\]
		Take point $\bar{x}$, $\bar{x}_i=1$, $i=1,\ldots,n$, where $f$ exhibits a kink and it is $f(\bar{x})=1$. Observe that none among the $2n$ (signed) coordinate  directions $\pm e_i$ is a descent one at $\bar{x}$ (it is in fact $f^\circ(\bar{x};-e_i)=0$ and $f^\circ(\bar{x};e_i)=1$, $i=1,\ldots,n$). Calculation of the $2n$ ratios $s_i$ as  in (\ref{incremental}), along the directions $e_i$ and $-e_i$ leads to $s_i=1$ and $s_i=0$ , respectively, for $i=1,\ldots, n$. It is easy to verify that, letting $p=n$ in Algorithm \ref{kmeansAlgo}, an optimal solution to problem (\ref{minmin}) is $\hat v_j=e_j$, $j=1,\ldots,n$. Finally, solving 
		\[\bar{d}=-\argmin_{v \in conv\{\hat v_j ,\ j=1,\ldots,n\}} \|v\|\]
		we obtain $\bar{d}=\displaystyle -\frac{e}{n} $, which is indeed a descent direction at $\bar{x}$.
	}
\end{example}

\section{The improved CS-DFN algorithm}\label{sec:improvedalgo}

This section is devoted to the definition of the improved version of algorithm CS-DFN which we call Fast-CS-DFN. The method is basically the CS-DFN Algorithm introduced in reference \cite{Fasano:14}, a derivative-free linesearch-type algorithm for the minimization of black-box (possibly) nonsmooth functions. It works by performing derivative-free linesearches along the coordinate directions and resorting to the use of a further search direction when the stepsizes used to explore the coordinate directions are sufficiently small. The rationale behind this choice is connected with the observation that the coordinate directions might not be descent directions near a non-stationary point of non-smoothness. In such situations, a richer set of directions must be used to (at least asymptotically) be able to improve the non-stationary point. The convergence anaysis of CS-DFN carried out in \cite{Fasano:14} hinges on the use of asymptotically dense sequences of search directions so that, at non-stationary points, for sufficiently large $k$ a direction of descent is used.

The algorithm that we propose, namely Fast-CS-DFN, is a modification of CS-DFN. The relevant differences between the two methods are:
\begin{enumerate}
	\item for the sake of simplicity, problem (\ref{prob}) is unconstrained; hence in Fast-CS-DFN no control to enforce feasibility with respect to the bound constraints is needed;
	
	\item after the deployment of the direction $d_k$, Fast-CS-DFN makes use of Algorithm \ref{DircoAlgo} to compute a direction that tries to exploits the information gathered during the optimization process to heuristically improve the last produced point.
\end{enumerate}

The Fast-CS-DFN Algorithm is reporten in Algorithm \ref{fast_csdfn}.

\begin{algorithm}[H]
	\caption{{\bf  Algorithm Fast-CS-DFN}}\label{fast_csdfn}
	
	\begin{algorithmic}[1]
		
		\par\medskip
		\State {\bf Input.} $\theta \in (0,1)$, $\eta > 0$, $x_0\in \Re^n$, $\tilde\alpha_0 > 0$, $\tilde\alpha_0^i > 0$, $d_0^i=e^i$, for $i=1,\ldots,n$, $G_0 = \emptyset$, 
		a sequence $\{d_k\}$ of search directions such that $\|d_k\|=1$, for all $k$.
		
		\For{$k=0,1,\dots$}
		
		\State Set $y_k^1=x_k$, $G^1_k = G_k$.
		
		\For{$i=1,2,\dots,n$}
		
		\State Compute $\alpha$, $d^i_{k+1}$, $G^{i+1}_k$ by the {\em Continuous Search}$(\tilde\alpha_k^i,y_k^i,d_k^i,G^i_k;\alpha, d^i_{k+1})$.
		
		
		\State {\bf If} $(\alpha = 0)$ {\bf then} set $\alpha_k^i = 0$ and $\tilde\alpha_{k+1}^i = \theta\tilde\alpha_k^i$
		
		\State {\bf else} set $\alpha_k^i = \alpha$ and $\tilde\alpha_{k+1}^i = \alpha$.
		
		
		\State Set $y_k^{i+1}=y_k^i+\alpha_k^id_{k+1}^i$.
		
		\EndFor
		
		\If{$\Big(\max_{i=1,\dots,n}\{\alpha_k^i,\tilde\alpha_k^i\} \leq \eta \Big)$}
		
		\State Compute $\alpha_k$, $\tilde d_k$ and $G_k^{n+2}$ by the {\em Continuous Search}$(\tilde\alpha_k,y_k^{n+1},d_k,G_k^{n+1};\alpha_k, \tilde d_k)$.

		\State {\bf If} $(\alpha_k = 0)$ {\bf then} $\tilde \alpha_{k+1}=\theta \tilde \alpha_k$ and $y_k^{n+2} = y_k^{n+1}$
		
		\State {\bf else} $\tilde\alpha_{k+1}= \alpha_k$ and $y_k^{n+2}=y_k^{n+1}+\alpha_k\tilde d_k$.
		
		
		\State \label{FCSDFN:step14} Build a symmetric matrix $B_k$
		\State \label{FCSDFN:step15}Compute  $\hat d_k^N$ by  {\em Direction computation}$(\tilde\alpha_k,y_k^{n+2},B_k,G_k^{n+2})$, i.e. Algorithm \ref{DircoAlgo}.
		
		\State \label{FCSDFN:step16}Compute $\check\alpha_k$, $\check d_k$ and $\bar G$ by the {\em Continuous Search}$(\tilde\alpha_k,y_k^{n+2},\hat d_k,G_k^{n+2};\check\alpha_k, \check d_k)$.

		\State \label{FCSDFN:step17}Set $x_{k+1} = y_k^{n+2}+\check\alpha_k \check d_k$
		
		\State \label{FCSDFN:step18}{\bf If} $\bar G = \emptyset$ {\bf then} set $G_{k+1} = \emptyset$ {\bf else} set $G_{k+1}= G_k^{n+2}$
		
		\Else
		
		\State set $\tilde\alpha_{k+1} = \tilde\alpha_k$, $\alpha_k = \check\alpha_k = 0$   and $y_k^{n+2}=y_k^{n+1}$, $G_k^{n+2} = G_k^{n+1}$.
		
		\State Set $x_{k+1} = y_k^{n+2}$, $G_{k+1} = G_k^{n+2}$		
		
		\EndIf
		

		\EndFor 
		
		\State {\bf Output.} The sequences $\{x_k\}$, $\{\alpha_k\}$, $\{\tilde\alpha_k\}$, $\{\alpha^i_k\}$ and $\{\tilde\alpha^i_k\}$, for $i=1,\dots,n$.
		
	\end{algorithmic}
\end{algorithm}

Some comments about Algorithm Fast-CS-DFN are in order.
\begin{enumerate}
	\item Fast-CS-DFN except for steps \ref{FCSDFN:step14}--\ref{FCSDFN:step18} and for the mechanism used to produce $G_{k+1}$ starting from $G_k$, exactly is the CS-DFN method as described in \cite{Fasano:14};
	
	\item the new direction $\hat d_k^N$ is used when the stepsizes $\alpha_k^i$ and $\tilde\alpha_k^i$, $i=1,\dots,n$, are sufficiently small and after the deployment of the direction $d_k$;
	
	\item the computation of the new direction $\hat d_k^N$ performed at step \ref{FCSDFN:step15} hinges (a) on the matrix $B_k$ and (b) on the set of couples $G_k^{n+2}$. 
	\begin{itemize}
		\item[(a)] To build $B_k$, we maintain a set of points $Y_k$ which is managed in just the same way as described in \cite{Fasano:14};
		
		\item[(b)]As for the set $G_k^{n+2}$, it stores information on the consecutive failures encountered up to the current point, i.e. in the deployment of the coordinate directions and the direction $d_k$. This set is emptied every time a non-null step is computed by the algorithm along any direction;
	\end{itemize} 
	
	\item the asymptotic convergence properties of Fast-CS-DFN are analogous to that of CS-DFN. The theoretical analysis follows quite easily from the results proved for CS-DFn in \cite{Fasano:14} when considering that the new iterate $x_{k+1}$ is such that $f(x_{k+1})\leq f(y_k^{n+2})$. 
\end{enumerate}

\begin{algorithm}[H]
	\caption{{\bf Continuous Search ($\tilde\alpha,y,p,G;\alpha,p^+$)}}\label{cont_search}
	
	\begin{algorithmic}[1]
		
		\par\medskip
		
		\State {\bf Data}: $\gamma > 0$, $\delta\in (0,1)$
		
		\State 
		Set $\alpha = \tilde\alpha$. 
		
		\State \label{CS:step1}
		{\bf If} $f(y+\alpha p)\le f(y)-\gamma\alpha^2$ {\bf
			then} set $p^+= \phantom{-}p$ and go to Step \ref{CS:step4}. 
		
		\State \label{CS:step2}
		{\bf If} $f(y-\alpha p)\le f(y)-\gamma\alpha^2$ {\bf then} set $p^+ = -p$ and go to Step
		\ref{CS:step4}.  
		
		\State \label{CS:step3}
		$G^+ \leftarrow G \cup \{(p,(f(y+\alpha p)-f(y))/\alpha),(-p,(f(y-\alpha p)-f(y))/\alpha)\} $
		\\ 
		Set $\alpha = 0$, {\bf return} $\alpha$, $p^+=p$ and $G^+$
		
		\State \label{CS:step4}
		Let $G^+\leftarrow\emptyset$ and $\beta={\alpha}/{\delta}$. 
		
		\State \label{CS:step5}
		{\bf If} $f(y+\beta p^+)>f(y)-\gamma\beta^2$  {\bf return} $\alpha, p^+$ and $G^+$ 
		
		\State \label{CS:step6}
		Set $\alpha=\beta$ and go to Step \ref{CS:step4}.
		
	\end{algorithmic}
	
\end{algorithm}


\section{Numerical results}\label{sec:results}

The proposed Fast-CS-DFN algorithm has been implemented in Python 3.9 and compared with CS-DFN \cite{Fasano:14} (available through the DFL library). The comparison has been carried out on a set of 47 nonsmooth problems. In the following subsections we briefly describe the test problems collection, the metrics adopted in the comparison and, finally, the obtained results.

\subsection{Test problems collection}
In Table \ref{tab1:prob} description of the test problems is reported. In particular, each table entry gives the problem name, the number $n$ of variables and the reference where the problem definition can be found.

\begin{table}[ht]
	{
		\begin{center}
			\begin{tabular}{|l|c|c|}\hline
				Problem name & $n$ & origin \\\hline\hline
				cb2 &   2 & \cite{lukvle:98} \\\hline
				crescent &   2 & \cite{karmitsa} \\\hline
				demymalo &   2 & \cite{lukvle:98} \\\hline
				davidon2 &   4 & \cite{lukvle:98} \\\hline
				kowalik &   4 & \cite{lukvle:98} \\\hline
				lukgamma &   4 & \cite{lukvle:98} \\\hline
				oet5 &   4 & \cite{lukvle:98} \\\hline
				oet6 &   4 & \cite{lukvle:98} \\\hline
				polak6 &   4 & \cite{lukvle:98} \\\hline
				colville1 &   5 & \cite{lukvle:98} \\\hline
				hs78 &   5 & \cite{lukvle:98} \\\hline
				lukexp &   5 & \cite{lukvle:98} \\\hline
				pbc1 &   5 & \cite{lukvle:98} \\\hline
				shor &   5 & \cite{lukvle:98} \\\hline
				elattar &   6 & \cite{lukvle:98} \\\hline
				evd61 &   6 & \cite{lukvle:98} \\\hline
				transformer &   6 & \cite{lukvle:98} \\\hline
				wong1 &   7 & \cite{lukvle:98} \\\hline
				lukfilter &   9 & \cite{lukvle:98} \\\hline
				gill &  10 & \cite{lukvle:98} \\\hline
				maxquad &  10 & \cite{lukvle:98} \\\hline
			\end{tabular}
			\begin{tabular}{|l|c|c|}\hline
				Problem name & $n$ & origin \\\hline\hline
				polak2 &  10 & \cite{lukvle:98} \\\hline
				wong2 &  10 & \cite{lukvle:98} \\\hline
				osborne2 &  11 & \cite{lukvle:98} \\\hline
				polak3 &  11 & \cite{lukvle:98} \\\hline
				steiner2 &  12 & \cite{lukvle:98} \\\hline
				shelldual &  15 & \cite{lukvle:98} \\\hline
				watson &  20 & \cite{lukvle:98} \\\hline
				wild1 & 20 & \cite{wild:dp} \\\hline
				wild2 & 20 & \cite{wild:dp} \\\hline
				wild3 & 20 & \cite{wild:dp} \\\hline
				wild11 & 20 & \cite{wild:dp} \\\hline
				wild15 & 20 & \cite{wild:dp} \\\hline
				wild16 & 20 & \cite{wild:dp} \\\hline
				wild19 & 20 & \cite{wild:dp} \\\hline
				wild20 & 20 & \cite{wild:dp} \\\hline
				wild21 & 20 & \cite{wild:dp} \\\hline
				wong3 &  20 & \cite{lukvle:98} \\\hline
				cb3 &  20,30,40 & \cite{karmitsa} \\\hline
				l1hilb &  20,30,40 & \cite{karmitsa} \\\hline
				maxq &  20,30,40 & \cite{karmitsa} \\\hline
				\multicolumn{3}{c}{ }\\
			\end{tabular}
		\end{center}
		\caption{Description of the test problems}\label{tab1:prob}
	}
\end{table}

\subsection{Metrics}
To compare our derivative-free algorithms we resort to the use of the well-known performance and data profiles (proposed in \cite{more:pp} and \cite{wild:dp}, respectively). In particular, let $\cal P$ be a set of problems and $\cal S$ a set of solvers used to tackle problems in $\cal P$. Let $\tau > 0$ be a required precision level and denote by $t_{ps}$ the performance index, that is the number of function evaluations required by solver $s\in\cal S$ to solve problem $p\in\cal P$. Problem $p$ is claimed to be solved when a point $x$ has been obtained such that the following criterion is satisfied
\[ 
f(x) \leq f_L + \tau(f(x_0)-f_L)
\]
where $f(x_0)$ is the initial function value and $f_L$ denotes  the best function value found by any solver on problem $p$ itself. Then, the performance ratio $r_{ps}$ is
\[
r_{ps} = \frac{t_{ps}}{\min_{i\in\cal S}\{t_{pi}\}}.
\]
Finally, the performance and data profiles of solver $s$ are so defined
\[
\rho_s(\alpha) = \frac{1}{|\cal P|}|\{p\in P:\ r_{ps} \leq\alpha\}|,\quad 
d_s(\kappa) = \frac{1}{|\cal P|}|\{p\in P:\ t_{ps}/(n_p+1) \leq\kappa\}|
\]
where $n_p$ is the number of variables of problem $p$. Particularly, the performance profile $\rho_s(\alpha)$ tells us the fraction of problems that solver $s$  solves with a number of function evaluation which is at most $\alpha$ times the number of function evaluations required by the best performing solver on that problem. On the other hand, the data profile $d_s(\kappa)$ indicates the fraction of problems solved by $s$ with a number of function evaluations which is at most equal to $\kappa(n_p+1)$, that is the number of function evaluations required to compute $\kappa$ simplex gradients.  

When using performance and data profiles for benchmarking derivative-free algorithms, it is quite usual to consider (at least) three different levels of precision (low, medium and high) corresponding to $\tau = 10^{-1},10^{-3},10^{-5}$, respectively. 

\subsection{Results} 
\begin{figure}[htb]
	\begin{center}
		\includegraphics[width=0.9\textwidth]{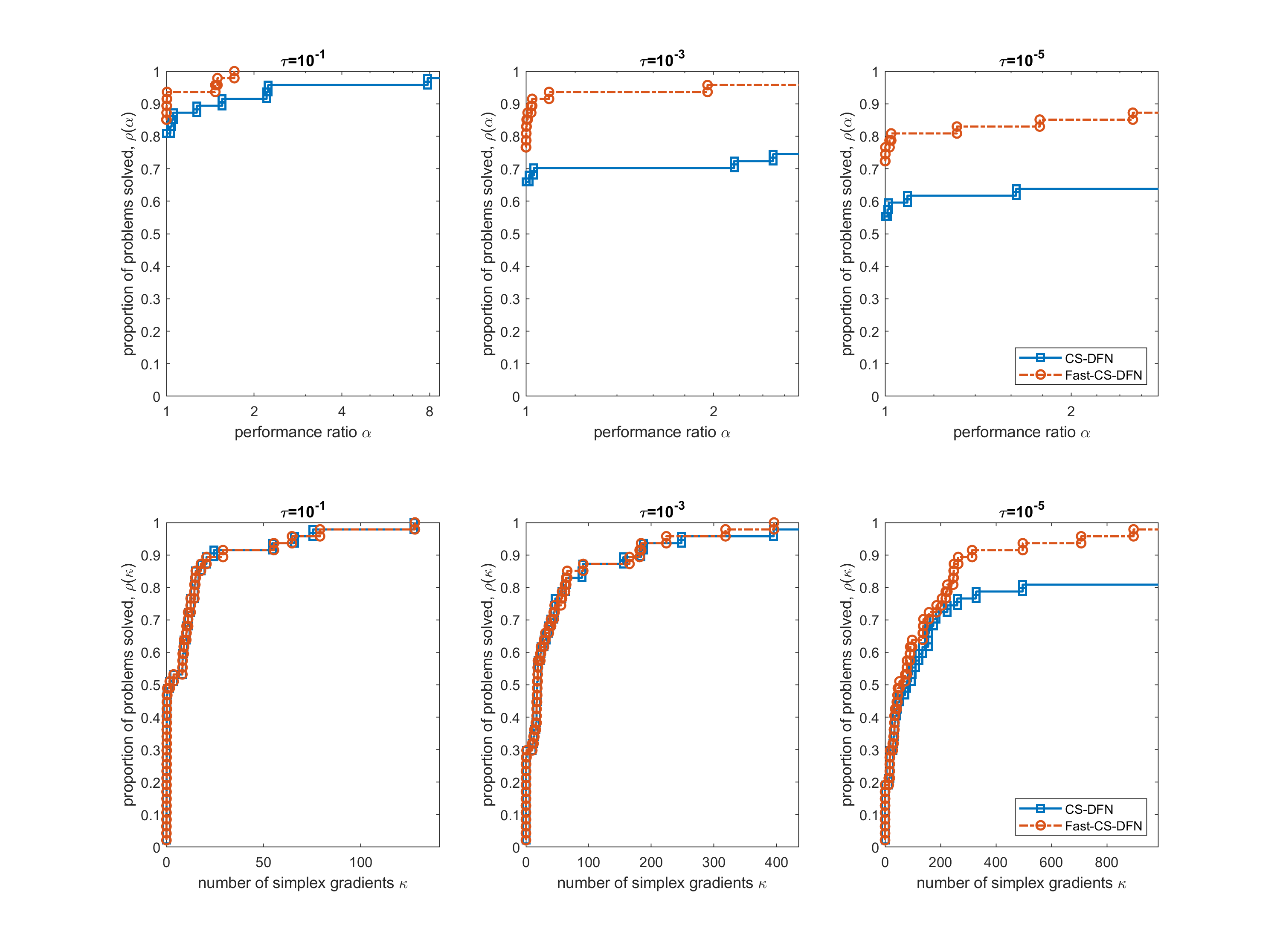}
	\end{center}
	\caption{Comparison of Fast-CS-DFN and CS-DFN}\label{fig1}
\end{figure}

Figure \ref{fig1} reports the results of the comparison by means of performance and data profiles between Fast-CS-DFn and CS-DFN. 

As we can see, the new algorithm Fast-CS-DFN is always more robust, namely it is able to solve the largest portion of problems within a given amount of computational effort. More in particular, from the performance profiles, we can also say that the new method is  invariably more efficient than the original one since the profile curves always have higher values for $\alpha = 1$.

\section{Conclusions}\label{sec:conclusion}
In the paper, we propose a strategy to compute  (possibly) good descent directions that can be further heuristically exploited within derivative-free algorithms for nonsmooth optimization. In fact, we show that the use of the proposed direction within the CS-DFN algorithm \cite{Fasano:14} improves the performances of the method. Numerical results on a set of nonsmooth optimization problems from the literature show the efficiency of the proposed direction computation strategy.

As a final remark, we point out that the proposed strategy could be embedded in virtually any optimization algorithm as an heuristic to try and produce improving points.   

\section*{Data availability statements}
The datasets generated during and/or analysed during the current study are available in the DFL repository, \texttt{http://www.iasi.cnr.it/$\sim$liuzzi/dfl} as package {\tt FASTDFN}.


\end{document}